\documentclass[10pt]{article}
\usepackage{amsfonts}
\usepackage{amssymb}
\usepackage{amsmath}
\usepackage{latexsym}
\usepackage{graphicx}

\usepackage{hyperref}
\usepackage{color}

\usepackage{multicol}
\usepackage{multirow}
\usepackage{colortbl}
\usepackage[bottom]{footmisc}
\usepackage{url}
\usepackage{datetime}

\usepackage{bm}

\usepackage[linesnumbered,ruled,vlined]{algorithm2e}
\usepackage{algorithmicx}
\usepackage{algpseudocode}
\usepackage{varwidth}

\usepackage{color}

\usepackage{geometry}

\newtheorem{thm}{Theorem}[section]
\newtheorem{lem}{Lemma}[section]

\newcommand{\be}{\begin{equation}}
\newcommand{\ee}{\end{equation}}
\newcommand{\ba}{\begin{array}}
\newcommand{\ea}{\end{array}}
\newcommand{\ben}{\begin{eqnarray}}
\newcommand{\een}{\end{eqnarray}}
\newcommand{\bn}{\begin{eqnarray*}}
\newcommand{\en}{\end{eqnarray*}}

\newcommand{\sfrac}[2]{\mathchoice
  {\kern0em\raise.5ex\hbox{\the\scriptfont0 #1}\kern-.15em/
   \kern-.15em\lower.25ex\hbox{\the\scriptfont0 #2}}
  {\kern0em\raise.5ex\hbox{\the\scriptfont0 #1}\kern-.15em/
   \kern-.15em\lower.25ex\hbox{\the\scriptfont0 #2}}
  {\kern0em\raise.5ex\hbox{\the\scriptscriptfont0 #1}\kern-.2em/
   \kern-.15em\lower.25ex\hbox{\the\scriptscriptfont0 #2}}
  {#1\!/#2}}

\renewcommand{\div}{{\rm div\,}}

\def\ub {{\bf u}}
\def\vb {{\bf v}}

\def\ub {{\bf u}}

\newcommand{\ds}{\displaystyle}

\newtheorem{my assumption}{Assumption}
\begin{document}
\centerline{\Large \bf Hybrid and Multiplicative Overlapping Schwarz Algorithms}
\vskip .2cm
\centerline{\Large \bf with Standard Coarse Spaces for}
\vskip .2cm
\centerline{{\Large \bf Mixed Linear Elasticity and Stokes Problems}}

\vskip .5cm
\renewcommand{\thefootnote}{\arabic{footnote}}
\centerline{{\bf Mingchao Cai}$\dagger$\footnote{Corresponding author. E-mail address:
cmchao2005@gmail.com.} and {\bf Luca F. Pavarino}$\S$\footnote{E-mail address: luca.pavarino@unimi.it. The work by this author was supported in part by M.I.U.R. grant (PRIN 201289A4LX\_002)}
}

\vskip .3cm \centerline{$\dagger$\it Department of Mathematics, Morgan State University,}
\centerline{\it 1700 E Cold Spring Ln, Baltimore, MD 21251, USA.}

\vskip .3cm \centerline{$\S$\it Dipartimento di Matematica, Universit\`a di Milano,}
\centerline{\it Via Saldini 50, 20133 Milano, Italy.}

\vskip .5cm \noindent\rule[2mm]{\textwidth}{.1pt}

\noindent {\bf Abstract.}
\vskip .5cm
The goal of this work is to construct and study hybrid and multiplicative two-level overlapping Schwarz algorithms with standard coarse spaces for the almost incompressible linear elasticity and Stokes systems, discretized by mixed finite and spectral element methods with discontinuous pressures. Two different approaches are considered to solve the resulting saddle point systems:
a)  a preconditioned conjugate gradient (PCG) method applied to the symmetric positive definite reformulation of the almost incompressible linear elasticity system obtained by eliminating the pressure unknowns;
b) a GMRES method with indefinite overlapping Schwarz preconditioner applied directly to the saddle point formulation
of both the elasticity and Stokes systems.
Condition number estimates and convergence properties of the proposed hybrid and multiplicative overlapping Schwarz algorithms are proven for the positive definite reformulation of almost incompressible elasticity.
These results are based on our previous study  \cite{cai2015overlapping} where only additive Schwarz preconditioners were considered for almost incompressible elasticity.
Extensive numerical experiments with both finite and spectral elements show that the proposed overlapping Schwarz preconditioners are scalable, quasi-optimal in the number of unknowns across individual subdomains and robust
with respect to discontinuities of the material parameters across subdomains interfaces.
The results indicate that the proposed preconditioners retain a good performance also when the quasi-monotonicity assumption, required by the available theory, does not hold.



\definecolor{kugray5}{RGB}{224,224,224}
\setlength{\marginparwidth}{0.75in}
\newcommand{\MarginPar}[1]{\marginpar{\vskip-\baselineskip\raggedright\tiny\sffamily\hrule\smallskip{\color{red}#1}\par\smallskip\hrule}}


\section{Introduction}

Finite and spectral element discretizations of the linear elasticity system suffer increasingly from locking effects and ill-conditioning, when the material approaches the incompressible limit, if only the displacement variables are used. One remedy to overcome the locking effect is using the following mixed form of the linear elasticity operator
\begin{equation}\label{elastic_operator}
\left[\begin{array}{cc}
-2 \mu \div {\bf D}(\cdot)  & \mbox{grad}  \\
-\div & -\frac{1}{\lambda} I
\end{array}\right].
\end{equation}
Here, $\lambda$ and $\mu$ are the Lam\'{e} constants, expressed as
\begin{equation}\label{Poisson}
\lambda = \frac{E \nu}
{(1+\nu)(1-2\nu)},
\quad  \mu=  \frac{E}
{2(1+\nu)}
\end{equation}
with $E$ being the modulus of elasticity (Young's modulus) and $\nu$ being the Poisson ratio of the elastic material. It is easy to see that the above mixed operator degenerates to the incompressible Stokes operator when the $(2, 2)$ block is zero. Finite and spectral element discretization of this mixed formulation lead to large saddle point systems whose iterative solution requires effective and efficient preconditioners. The goal of this paper is to construct and study hybrid and multiplicative two-level overlapping Schwarz preconditioners with standard coarse spaces for the mixed discretization of almost incompressible linear elasticity and Stokes systems. Earlier works on overlapping Schwarz methods for linear elasticity have focused on the compressible case in which the Poisson ratio $\nu$ is bounded away from $1/2$. The related recent developments are as follows: some nonoverlapping domain decomposition algorithms for mixed elasticity and Stokes systems have studied Wirebasket and Balancing Neumann-Neumann methods, see \cite{beir2006positive, goldfeld2003balancing, pavarino2000iterative, pavarino2002balancing}, and  FETI-DP and BDDC methods for the incompressible limit, see \cite{kim2012two, klawonn2007dual, li2005dual, li2005bddc, pavarino2010bddc, tu2013unified}; more recent applications include fluid-structure interaction \cite{barker2010two}, computational fluid dynamics \cite{hwang2007class, kim2012two}, and isogeometric analysis \cite{beir2013isogeometric}.

We will consider two different approaches to solve the resulting saddle point systems:
a)  a preconditioned conjugate gradient (PCG) method applied to the symmetric positive definite reformulation of the almost incompressible linear elasticity system obtained by eliminating the pressure unknowns element by element;
b) a GMRES method with indefinite overlapping Schwarz preconditioner applied directly to the saddle point formulation
of both the elasticity and Stokes systems. In both approaches, our main interest is the almost incompressible case and the incompressible Stokes limit.

Earlier works proposed some additive and hybrid two-level overlapping Schwarz algorithms in which exotic coarse spaces are used in the designing of the coarse grid problem, see \cite{dohrmann2009overlapping, dohrmann2010hybrid}. The tools and theoretical analysis developed in these papers
yield a condition number bound on the preconditioned additive Schwarz operator which is cubic (in contrast, it is linear for the compressible case) in the relative overlap size and which grows logarithmically with the number of elements across individual subdomains. In our previous work \cite{cai2015overlapping}, we have proposed and analyzed a two-level additive overlapping Schwarz algorithm in which standard coarse spaces are used. The corresponding condition number estimate obtained in \cite{cai2015overlapping} is cubic as in \cite{dohrmann2009overlapping}. However, the so-called {\it quasi-monotonicity} (see Section 3) assumption on the coefficient distribution is required when developing the theory. Inspired by our previous work \cite{cai2015overlapping} , in this paper we develop hybrid and multiplicative overlapping Schwarz algorithms using standard coarse spaces. By using some classical results \cite{mandel1994hybrid, toselli2005domain} and the newly developed theory from \cite{cai2015overlapping}, we analyze the proposed hybrid and multiplicative Schwarz algorithms. Both theoretical analysis and numerical experiments show that all our two-level algorithms are robust with respect to the number of subdomains, their diameters and mesh sizes and possible discontinuities of the material parameters across the subdomain interfaces.
In comparison with the exotic coarse space two-level methods developed in \cite{dohrmann2009overlapping, dohrmann2010hybrid}, the numerical results for our method also seem to indicate a faster convergence rate without the logarithmic factor $(1+\log(H/h))$ when  $H/\delta$ is fixed, required by the theoretical convergence estimates.
Furthermore, our numerical experiments show that the proposed preconditioners seem to perform very well also for some cases in which the {\it quasi-monotonicity} assumption does not hold, indicating the possibility of improving the condition number estimates or weakening the {\it quasi-monotonicity} assumption. We believe that our experimental investigation will provide valuable clues for future theoretical studies.



The second approach based on indefinite overlapping Schwarz preconditioners has been considered in the earlier works \cite{fischer1997overlapping, klawonn1998overlapping, klawonn2000comparison, pavarino1998preconditioned, pavarino2000indefinite}. These preconditioners are constructed directly from the indefinite saddle point systems and are based on both local and coarse saddle point problems. In this work, we discuss how to select the associated functional space of the subdomain problems. It is found that the pressure function spaces in the construction of the subdomain solvers are too large if we use the traditional zero mean value pressure spaces,
resulting in small local inf-sup constants
and preconditioners with slower convergence rates.  Numerical experiments show that in addition to the traditional zero mean value  constraint for the pressure subspaces, one should also impose zero Dirichlet boundary condition on the extended layers degrees of freedom of each subdomain. Although the supporting theory for these indefinite preconditioners is still lacking,  the proposed algorithms are scalable, efficient and robust with respect to the material incompressibility and possible discontinuities of the material parameters $\lambda, \mu$ across subdomain boundaries.

The rest of this paper is organized as follows. In Section 2, we introduce the mixed finite and spectral element approximation of the linear elasticity and Stokes problems. In Section 3, we study the symmetric positive definite reformulation of the linear elasticity problem and we introduce additive, hybrid and multiplicative overlapping Schwarz algorithms for it. In Section 4, indefinite overlapping Schwarz preconditioners are designed for the saddle point formulation. In Section 5, extensive numerical experiments based on the ${\bf Q}_2^h-P_1^h$ elements and the $Q_{n}-Q_{n-2}$ spectral elements are presented. The scalability of the two-level algorithms are numerically verified. We also check the dependence of the condition number on the relative overlapping size, the number of elements in each individual subdomain and the preconditioner robustness with respect to discontinuous coefficients  and checkerboard tests.

\section{Mixed formulation of linear elasticity and Stokes problems}

We consider a domain $\Omega \subset R^d, d=2,3$, decomposed into $N$ nonoverlapping subdomains $\Omega_i$ of diameter $H_i$, forming a coarse finite element partition $\tau_H$ of $\Omega$,
\begin{equation}
\overline{\Omega} = \bigcup_{i=1}^{N} \overline{\Omega}_i.
\label{H_subdomains}
\end{equation}
Here, $H = \max_i H_i$ is the characteristic diameter of the subdomains. To simplify our discussion, we assume that all boundary conditions are of Dirichlet type and then $\partial \Omega_D = \partial \Omega$. We further assume that the solution vanishes on $\partial \Omega$. We denote the interface of the domain decomposition (\ref{H_subdomains}) as
$$
\Gamma = \left(\bigcup_{i=1}^N \partial \Omega_i \right)
 \setminus \partial \Omega.
$$
A fine triangulation ${\mathcal T}_h$ of $\Omega$ is obtained by partitioning each subdomain into many shape-regular finite elements. We will assume that the nodes match across the interface between the subdomains. We will denote by $H/h$ the number of elements on each subdomain side (without overlap). For each $\Omega_i$, we obtain a larger subdomain $\Omega_i^\prime$ by adding layers of elements around its boundary. We will denote the minimal thickness of $\Omega_i^\prime \setminus \Omega_i$ by $\delta_i$. Then, with the above introduction, $\ds \overline{\Omega}=\bigcup_{i=1}^N \overline\Omega_i^\prime$ is the overlapping partition of $\Omega$.

We consider the standard displacement and pressure spaces
\[ \bm{V} := \{\bm{v} \in
H^1(\Omega)^d :\bm{v}|_{\partial\Omega_D} = 0\},
\quad Q:= L^2_0(\Omega) := \{p \in L^2(\Omega),  \int_\Omega p dx =0\},
\]
and the following mixed formulation of the linear elasticity problem (see e.g. \cite[Chapter 1]{boffi2013mixed}):\\
find $(\bm{u},p) \in \bm{V}\times Q$ such that
\begin{equation}
\left\{\begin{array}{cccccl}
a(\bm{u}, \bm{v})  & + &
b(\bm{v}, p) & = &
 \langle{\bf F},\bm{v}\rangle &
\forall \bm{v} \in \bm{V}, \\
 & & & & & \\
b(\bm{u}, q)  & + &
c(p,q) & = & 0 & \forall q \in Q , \\
\end{array}\right.
\label{mixed_el}
\end{equation}
with bilinear forms defined by assembling local contributions from the different subdomains
\begin{eqnarray}
a(\bm{u}, \bm{v}) & = & \sum_{i=1}^{N}\mu_i a_i(\bm{u},\bm{v}) :=
2\sum_{i=1}^{N}\mu_i\int_{\Omega_i}{\bf D}(\bm{u}):{\bf D}(\bm{v})\ dx ,
\label{bilinear_forms_a}\\
c(p, q) & = & \sum_{i=1}^{N}\frac{1}{\lambda_i} c_i(p,q):=
\sum_{i=1}^{N}\frac{1}{\lambda_i}\int_{\Omega_i}p\ q\ dx ,
\label{bilinear_forms_c}\\
b(\bm{v}, q) & = & \sum_{i=1}^{N} b_i(\bm{v},q) :=
-\sum_{i=1}^{N}\int_{\Omega_i}{\rm div}\bm{v}\ q\ dx .
\label{bilinear_forms_b}
\end{eqnarray}
Here,
$$
\ds
{\bf D}({\ub}) :{\bf D}({\vb}) =
\sum_{i=1}^d
\sum_{j=1}^d
{\bf D}_{ij}(\ub) {\bf D}_{ij}(\vb),
\quad \mbox{with} \quad
{\bf D}_{ij}(\ub) = \frac{1}{2} \left(\frac{\partial u_i}{\partial x_j}
+ \frac{\partial u_j}{\partial x_i} \right),
$$
${\bf F}$ represents the applied forces and for simplicity we assume constant  Lam\'{e} parameters in each subdomain $\Omega_i$, i.e., $\mu = \mu_i$ and $\lambda = \lambda_i$ in $\Omega_i$. These parameters can be expressed in terms of the local Poisson ratio $\nu_i$ and Young's
modulus $E_i$ as those in (\ref{Poisson}).

If all $\lambda_i \rightarrow \infty$, we have  $c(p,q)=0$ and we obtain the limiting problem for incompressible linear elasticity or the classical Stokes system for an incompressible fluid with pure Dirichlet boundary condition.
Generalized Stokes problems originating from stabilization techniques or penalty formulations can also be written as in
the model problem (\ref{mixed_el}).

{\bf ${\bf Q}_2^h$ discontinuous $P_1^h$ mixed finite elements.}
We now introduce the ${\bf Q}_2^h$ discontinuous $P_1^h$ finite element approximation of (\ref{mixed_el}). Let $K$ be a rectangle of ${\mathcal T}_h$, the conforming ${\bf Q}_2^h-P_1^h$ finite element space for the displacement is
\bn
\bm{V}^h :=\big\{\bm{v} \in \bm{V} \big|\ \bm{v}|_K \in Q_2(K)^d,\quad d=2,3, ~\forall K\in{\cal T}_{h} \big\}.
\en
and the pressure space
\bn
Q^h :=\big\{q\in L^2_0(\Omega) \big|\ q|_K \in P_1(K),\ \forall K\in{\cal T}_{h} \big\}.
\en
In 2D, for each element, there are 9 local degrees of freedom for the displacement, located at the vertices, midpoints of the edges and in the center of the quadrilateral; and there are 3 degrees of freedom (with the function value and the two partial derivatives) located at the center of the quadrilateral. In 3D, there are 27 degrees of freedom for displacement located at the cell vertices, the midpoints of the cell edges, the middle of the six faces of the cube, the center of the element, and there are 4 pressure degrees of freedom located at the center of the element. The approximation property of the ${\bf Q}_2^h-P_1^h$ elements can be found in \cite{boffi2002quadrilateral, matthies2002inf}: the $H^1-$ norm of the displacement approximation is $\mathcal{O}(h^2)$, while the pressure approximation order depends on whether one use the mapped version or unmapped version of the  ${\bf Q}_2^h-P_1^h$ element. For the unmapped version, the $L^2-$ norm error is $\mathcal{O}(h^2)$, but for the mapped version, the $L^2-$ error order is $\mathcal{O}(h)$.

For the underlying finite element discretization, there exists a positive constant $\beta$, such that the following uniform inf-sup condition holds.
\begin{equation}
\sup_{\bm{v} \in \bm{V}^h}
\frac{b_i(\bm{v},q)}{a_i(\bm{v},\bm{v})^{1/2}}  \geq
\beta c_i(q,q)^{1/2} \ \  \forall q\in Q^h, 
\ \ \beta > 0  \nonumber.
\end{equation}
For a proof of this result, we refer to \cite{boffi2002quadrilateral, matthies2002inf} and \cite[pp.\ 156--158]{girault1986finite}.

{\bf ${\bf Q}_n$ discontinuous $Q_{n-2}$ mixed spectral elements.}
We also consider mixed spectral element discretizations, see e.g. \cite{BM1997, CHQZ2007, le1997non}, where
the space of displacements $\bm{V}$ is discretized, component by component, by
continuous, piecewise tensor product polynomials of degree $n\geq 2$:
\[
\bm{V}^h := \{\bm{v}\in \bm{V} : {v}|_{K_{i}}\circ \phi_i
\in Q_{n}(K_{\rm ref})^d,  \ d=2,3, \ n \ge 2, \forall K_{i}\in {\cal T}_{h}\} ,
\]
where each element  is an affine image of the reference cube $K_i = \phi_i(K_{\rm ref})$
with an affine map $\phi_i$.

The pressure space $Q$ is discretized by discontinuous, piecewise tensor product
polynomials of degree $n-2$:
\[
Q^h := \{q \in L^2_0(\Omega) : q|_{K_i}\circ \phi_i \in Q_{n-2}(K_{\rm ref}),
 \ n \ge 2, \forall K_{i}\in {\cal T}_{h}\} .
\]

We use Gauss-Lobatto-Legendre (GLL(n)) quadrature, which also allows for the construction of very convenient nodal tensor-product bases for $\bm{V}^h$ and $Q^h$, using for the latter only the interior GLL nodes of each element. We denote by $\{\xi_i\}_{i=0}^n$ the set of GLL(n) points of $[-1,1]$, by $\sigma_i$ the quadrature weight associated with $\xi_i,$ and by $l_i(x)$ the Lagrange
interpolating polynomial of degree $n$ that vanishes at all the GLL(n) nodes except at $\xi_i$, where it equals $1$. Each  element
of $Q_n(T_{\rm ref})$ is expanded in this GLL(n) basis, and any $L^2-$inner product of two scalar components $u$ and $v$ is
replaced by
\begin{equation}
(u,v)_{n,\Omega} =  \sum_{s=1}^{N_e}\sum_{i,j,k=0}^n (u\circ
\phi_s)(\xi_i,\xi_j,\xi_k) (v\circ \phi_s)(\xi_i,\xi_j,\xi_k)
|J_s|\sigma_i\sigma_j\sigma_k \,, \nonumber 
\end{equation}
where $|J_s|$ is the determinant of the Jacobian of $\phi_s$. Similarly, a very convenient basis for $Q_{n-2}$ consists of
the tensor-product Lagrangian nodal basis functions associated with the internal GLL(n) nodes, i.e., the endpoints $-1$ and $+1$ are excluded. The mass matrix based on these basis elements and GLL(n) quadrature is then diagonal for the displacement field but not for the pressure field.

The $Q_n-Q_{n-2}$ method satisfies a nonuniform inf-sup condition
\begin{equation}
\sup_{\bm{v}\in \bm{V}^h} \frac{({\rm div}\bm{v},q)}{\|\bm{v}\|_{H^1}}\geq \beta_n \|q\|_{L^2} \ \ \ \  \forall q\in U^h\, ,
\label{inf-supn}
\end{equation}
where $\beta_n \geq C n^{-1}$ and the constant $C > 0$ is independent of $n$. It is also known that $\beta_n$ decays slower for small $n$ than indicated by the theoretical bound; for example, $\beta_n\geq 0.43$ for $n\leq 16$ according to Maday et al. \cite{maday1993analysis}. (We note however, that the inf-sup coefficient will decrease with an increase in the aspect ratio of
a domain, see Dobrowolski \cite{dobro2003lbb}.) An alternative mixed spectral element method, with a bound on the inf-sup
constant which does not depend on $n$, is provided by the $Q_n-P_{n-1}$ method; see Bernardi and Maday \cite{bernardi1999uniform};

{\bf The discrete saddle point system.}
The discrete system obtained from the mixed finite or spectral elements introduced above is assembled from the saddle point operators of the subdomains $\Omega_i:$
\begin{equation}
\left[
\begin{array}{cc}
\mu_i A_i & B_i^{T} \\
B_i & - \frac{1}{\lambda_i}~C_i \\
\end{array}\right],
\label{local_saddle_pt}
\end{equation}
where $\mu_iA_i, B_i$, and $1/{\lambda_i}~C_i$ are the operators associated with the local bilinear forms $\mu_ia_i(\cdot,\cdot),\ b_i(\cdot,\cdot)$, and $1/{\lambda_i}c_i(\cdot,\cdot)$ defined in (\ref{bilinear_forms_a}), (\ref{bilinear_forms_c}), and (\ref{bilinear_forms_b}), respectively. If $c_i=0$, then we obtain the symmetric indefinite linear system for the Stokes equations.

\section{Overlapping Schwarz algorithms for the symmetric positive definite reformulation}


Since we are using discontinuous pressures, all pressure degrees of freedom can be eliminated, element  by element, to obtain reduced positive definite operators
\begin{equation}
\bar{A}_i := \mu_i A_i + \lambda_i B_i^{T}
C_i^{-1}B_i, \nonumber
\label{pos_def_Ai}
\end{equation}
that can be subassembled into a global positive definite operator $\bar{A}$. In case of constant coefficients $\mu, \lambda$, we have
\begin{equation}
\bar{A} := \mu A + \lambda B^{T}
C^{-1}B.
\label{reform_A}
\end{equation}
The bilinear form associated to $\bar{A}$, the operator of the positive definite reformulation, is denoted as $\bar{a}(\cdot,\cdot)$.

We now introduce the decomposition into local and coarse spaces. The coarse space on the coarse subdomain mesh $\tau_H$ is denoted by
\[
\bm{V}_0 = \bm{V}^H := \{\bm{v}\in \bm{V} : \bm{v}|_{\Omega_i} \in (Q_2(\Omega_i))^d~~ \forall \Omega_i \in \tau_H \}.
\]
The local problems are defined on the extended subdomains $\Omega_i^\prime$. To each of the $\Omega_i^\prime$, we associate a local space
$$
\bm{V}_i =  \bm{V}^h(\Omega_i^\prime)
\cap \bm{H}^1_0(\Omega_i^\prime)
$$
and a bilinear form $\bar{a}_i^\prime(\bm{u}_i,\bm{v}_i) := \bar{a}(R_i^T\bm{u}_i,R_i^T \bm{v}_i)$, where $R_i^T: \bm{V}_i \rightarrow  \bm{V}^h,$ simply extends any element of $\bm{V}_i$ by zero outside $\Omega_i^\prime$. Then, as we will only consider algorithms for which the local problems are solved exactly, we find that the local operators are
$$
\bar{A}_i^{\prime} = R_i \bar{A} R_i^T, \quad i=1, ..., N.
$$

Given the local and coarse embedding operators ${R}_{i}^T:\bm{V}_{i}\rightarrow \bm{V}^h$, $i=1,...,N$, and ${R}_0^T:\bm{V}_{0}\rightarrow \bm{V}^h$, the discrete space $\bm{V}^h$ can be decomposed into coarse and local spaces as
\[
\bm{V}^h = {R}_0^T \bm{V}_0   + \sum_{i}{R}_{i}^T\bm{V}_{i}.
\]
We denote the local and coarse (for $i=0$) projections as $P_{i}:={R}_{i}^T\bar{P}_{i} : \bm{V}^h \rightarrow {R}_{i}^T \bm{V}_i \subseteq \bm{V}^{h}$ with $\bar{P}_{i}: \bm{V}^h \rightarrow \bm{V}_i$ defined by
$$
\bar{a}_i^{\prime}(\bar{P}_{i} \ub,\vb_i)= \bar{a}(\ub,{R}_{i}^T \vb_i), \quad\forall \vb_i \in \bm{V}_{i}, \ \ \ i=0,\cdots,N.
$$
Note that we are using exact solvers for all the subspaces, we find that ${P}_{i}$ are all projections; cf.\ \cite[Section 2.2]{toselli2005domain}.

{\bf Additive Schwarz algorithm (OAS).}
The two-level Overlapping Additive Schwarz (OAS) operator is then given by
\begin{equation}\label{ADDSCH}
{P}_{OAS}:={P}_0+\sum_{i=1}^N{P}_{i}.
\end{equation}
We denote the operator form of the two-level additive Schwarz preconditioner as $B_{OAS}$, and we have ${P}_{OAS}=B_{OAS} \bar{A}$, where
\begin{equation}
B_{OAS}=R_0^T\,\bar{A}_0^{-1}\,R_0 + \sum_{i=1}^NR_{i}^T\,\bar{A}_i^{\prime{-1}}\,R_{i}.
\label{B_OAS}
\end{equation}
Here, $\bar{A}_{0} = R_0\bar{A} R_0^T$ is the coarse problem operator.

{\bf Hybrid Schwarz algorithm (OHS).}
Our hybrid Schwarz method is a variant of the additive Schwarz method. It is additive with respect to the local components and multiplicative with respect to the levels, see \cite{mandel1994hybrid,toselli2005domain}. More clearly, the action of the hybrid overlapping Schwarz algorithm on a given vector $\bm{r} \in \bm{V}^h$ is as follows.
\begin{itemize}
\item[]
1. Compute the coarse grid approximation: find $\bar{\bm{u}}_0 \in \bm{V}_0$, such that
$$
\bar{a}(\bar{\bm{u}}_0, \bm{v}_0) = <\bm{r}, \bm{v}_0>, \quad \forall \bm{v}_0 \in  \bm{V}_0.
$$
\item[]
2. Loop over subdomains (in parallel): find $\bm{u}_i \in \bm{V}_i, i=1, ..., N$, such that
$$
\bar{a} (\bm{u}_i, \bm{v}_i) = < \bm{r}, \bm{v}_i> - \bar{a} (\bar{\bm{u}}_0, \bm{v}_i), \quad \forall \bm{v}_i \in  \bm{V}_i,
$$
and define the sum $\bm{w}:= \bm{u}_1 + ... + \bm{u}_N$.
\item[]
3. Find $\tilde{\bm{u}}_0 \in \bm{V}_0$ such that
$$
\bar{a} (\bm{w}-\tilde{\bm{u}}_0, \bm{v}_0) = <\bm{r}, \bm{v}_0>, \quad \forall \bm{v}_0 \in  \bm{V}_0,
$$
and return the solution, $\bm{u} = \bm{w} - \tilde{\bm{u}}_0$.
\end{itemize}

We denote $\bm{u}=\bar{A}^{-1} \bm{r}$. As exact solvers are employed for both the coarse problem and the local subproblems, we firstly note that $\bar{\bm{u}}_0 = P_0 \bm{u}$. From the second step, we have
$$
\bm{u}_i= P_i(\bm{u}-\bar{\bm{u}}_0)=P_i(I-P_0) \bm{u},
$$
and $\bm{w}= \sum_{i=1}^N P_i (I-P_0) \bm{u}$. Then, from the third step, $\tilde{\bm{u}}_0= -P_0 (\bm{u}-\bm{w})$, and finally we have
$$
\bm{u}= \bm{w} + P_0 (\bm{u}-\bm{w}) = (I-P_0) \bm{w}+ P_0 \bm{u} = (I-P_0)\sum_{i=1}^N P_i (I-P_0)\bm{u} + P_0 \bm{u}.
$$
In conclusion, the overlapping hybrid Schwarz (OHS) operator is
\begin{equation}\label{ohs_preconditioner}
P_{OHS}=P_0 + (I-P_0)\sum_{i=1}^N P_i(I-P_0).
\end{equation}

{\bf Multiplicative Schwarz algorithm (OMS).}
As in \cite[Section 2.2]{toselli2005domain}, for a given vector $\bm{r} \in \bm{V}^h$, by sequentially finding approximation in each subdomain and the coarse level, the overlapping multiplicative Schwarz algorithm reads as follows.
\begin{itemize}
\item[]
1. Set ${\bm{u}}=\bm{0}$.
\item[]
2. For i=0, ..., N, find $\bm{u}_i \in \bm{V}_i$ such that
$$
\begin{array}{l}
\bar{a}(\bm{u}_i, \bm{v}_i) = <\bm{r}, \bm{v}_i>- \bar{a} (\bm{u}, \bm{v}_i), \quad \forall \bm{v}_i \in  \bm{V}_i; \\
\bm{u} \leftarrow \bm{u} + \bm{u}_i.
\end{array}
$$
\item[]
3. Return $\bm{u}$.
\end{itemize}

As a result, the overlapping multiplicative Schwarz operator is
$$
P_{OMS}=I - (I-P_N)(I-P_{N-1})\cdots(I-P_0):= I - E_{OMS},
$$
where $E_{OMS}$ is the error propagation operator.
Similarly to $B_{OAS}$, we can define the hybrid and multiplicative Schwarz preconditioners ${B}_{OHS}$ and ${B}_{OMS}$ (cf. \cite[Section 2.2]{toselli2005domain}).

\subsection{Condition number estimates for the symmetric positive definite reformulation}

For completeness, we briefly recall the assumptions and the theoretical analysis of the two-level overlapping additive Schwarz algorithm. The following assumption was introduced in \cite[Section 5]{dryja1996multilevel} in work on multi-level Schwarz algorithms on scalar elliptic problems; the role of the coefficient function for these elliptic problems will be played by the set of Lam\'{e} parameters $\{\mu_i\}$ for the case at hand.

\begin{figure}
\centering
\vspace{-0.5cm}
\includegraphics[height=5.2cm]{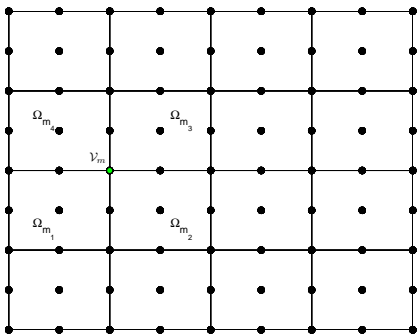}
\caption{An illustration of $\widehat{\Omega}_{m}$ in 2D with ${\bf Q}_2^H-P_1^H$ discretization.}
\label{trian}
\end{figure}

Let $\mathcal{V}_m$  be a vertex of the coarse grid $\tau_H$, with $m$ being the vertex global index. We denote $\widehat{\Omega}_{m}$ as the union of the elements in $\tau_H$ that share $\mathcal{V}_m$. For $\widehat{\Omega}_{m}$, let $N(m)$ be the total number of subdomains in $\widehat{\Omega}_{m}$ (cf. Figure \ref{trian} for an illustration).

\begin{my assumption} \label{ass1}
For each $\widehat{\Omega}_{m}$, order its subdomains such that $\mu_1=\mbox{max} \{ \mu_i \}$. We say that a distribution of the $\{ \mu_i \}_{i=1}^{N(m)}$ is quasi-monotone in $\widehat{\Omega}_{m}$ if for every $i$, there exists a sequence $\{i_j\}_{j=1}^{S}$, with
$$
\mu_{m_i} = \mu_{i_{S}}  \le, ..., \le \mu_{i_{j+1}} \le \mu_{i_{j}} \le, ..., \le  \mu_{i_1} = \mu_1,
$$
where the subdomains $\Omega_{i_j}$ and $\Omega_{i_{j+1}}$ have a face in common. If the vertex $\mathcal{V}_m$ belongs to $\partial \Omega$, then we additionally assume that $\partial \Omega_{1} \cap \partial \Omega$ contains a face for which $\mathcal{V}_m$ is a vertex.

A distribution $\{\mu_i\}$ on $\Omega$ is quasi-monotone with respect to the coarse triangulation $\tau_H$ if it is quasi-monotone for each $\widehat{\Omega}_{m}$.
\end{my assumption}

In the above assumption, we actually require that for any subdomain $\Omega_i$, there should exist a path, passing exclusively through subdomain faces of this set of subdomains, from $\Omega_i$ to the subdomain with maximum $\mu$ such that the values of the Lam\'{e} parameter $\mu_j$ are monotonically increasing along the path. For example, in Figure \ref{trian}, from the subdomain $\Omega_{m_3}$ to the subdomain $\Omega_{m_1}$, the path can be  clockwise or counter-clockwise.  The above assumption was relaxed in \cite{cai2015overlapping} by allowing a modest decrease when passing from one subdomain to the next along the path. In this work, we adopt the following assumption.
\begin{my assumption} \label{ass2}
Consider the same paths as in the previous assumption but relax the condition by assuming that for any subdomain $\Omega_j$ along the path from $\Omega_i$ to the subdomain with maximum $\mu$, the Lam\'{e} parameter satisfies $\mu_j \geq \mu_i$.
\end{my assumption}


For a counter example that does not satisfies either Assumption 1 or Assumption 2, we refer to Figure \ref{trian} and define the coefficient $\mu$ in $\Omega_{m_4}$ to $\Omega_{m_1}$ to be $(400, 600, 300, 700)$. Then, we see that $\Omega_{m_1}$ has the largest $\mu$. However, there is no path from $\Omega_{m_3}$ to $\Omega_{m_1}$ such that the $\mu$ values are monotonically increasing. For 2D cases with regular rectangular triangulations, Assumption 1 is actually equivalent to Assumption 2. For 3D cases,  Assumption 1 and Assumption 2 can be differentiated by considering the following two examples:
i) if we define the $\mu$ values along a path to be $(300, 400, 600, 700)$, then this distribution satisfies both Assumption 1 and Assumption 2;
ii) if we define the $\mu$ values along the same path to be $(300, 600, 400, 700)$, then this distribution satisfies Assumption 2 but not Assumption 1. 

As pointed out in \cite{cai2015overlapping}, these assumptions are made so that the union $\widehat{\Omega}_m$ of the subdomains and subdomain faces for all such paths associated with all the subdomain vertices of a single subdomain form a domain for which Poincar\'{e}'s and Korn's inequalities can still be used.

\begin{thm}\label{th:est_AS2level}
Let the set of Lam\'{e} parameters $\{\mu_i\}$ be quasi-monotone in the sense of Assumption \ref{ass2}. Then the condition number of the additive Schwarz operator $P_{OAS}$ satisfies
\begin{equation}\label{oas-cond}
\kappa(P_{OAS}) \leq C(H/\delta)^3(1+\log(H/\delta))(1+\log(H/h)),
\end{equation}
where $C$ is a constant independent of the number of subdomains, their diameters and mesh sizes, and which depends only on the number of colors required for the overlapping subdomains and the shape regularity of the elements and subdomains.
\end{thm}

Here, $H/\delta:= \max_i H_i/\delta_i$ and $H/h:= \max_i H_i/h_i$ where $H_i,h_i$ denote the
diameter and mesh size of $\Omega_i$ and $\delta_i$ measures the overlap of $\Omega_i^\prime$ and its next neighbors. For the analysis of the two-level overlapping additive Schwarz method, the main challenge is to develop a bound on the parameter $C_0^2$ as in \cite[Assumption 2.2]{toselli2005domain}:
\begin{equation}  \label{C-zero}
   \bar{a}(R_0^T\bm{u}_0,R_0^T\bm{u}_0) +
\sum_{i=1}^N \bar{a}_i^\prime(R_i^T\bm{u}_i,R_i^T\bm{u}_i)
\leq C_0^2 \bar{a}(\bm{u},\bm{u}),
\end{equation}
for any decomposition $\bm{u} = \sum_{i=0}^N R_i^T\bm{u}_i$. The detailed proof of this theorem can be found in \cite{cai2015overlapping}. We recall here that $C_0^2$ is defined as the right hand side of (\ref{oas-cond}) and $C_0^{-2}$ provided a lower bound for the eigenvalues of $P_{OAS}$.

The following lemma shows that the hybrid method yields improved convergence over the additive Schwarz method \cite{mandel1994hybrid}.
\begin{lem}\label{th:ohs_comp_oas}
The extreme eigenvalues of $P_{OHS}$ and $P_{OAS}$ satisfy
$$
\lambda_{min}(P_{OHS}) \ge \lambda_{min}(P_{OAS}), \quad \lambda_{max}(P_{OHS}) \le \lambda_{max}(P_{OAS}).
$$
\end{lem}

{\it Proof}. Let us study the following Rayleigh quotient
\begin{equation}
\label{Rayleigh_OHS}
\frac{\bar{a}(P_{OHS}\bm{u}, \bm{u})}{\bar{a}(\bm{u}, \bm{u})}=\frac{\bar{a}(P_0\bm{u}, P_0\bm{u})+\bar{a}(P_{OAS}(I-P_0)\bm{u}, (I-P_0)\bm{u})}{\bar{a}(P_0\bm{u}, P_0\bm{u})+\bar{a}((I-P_0)\bm{u}, (I-P_0)\bm{u})}
\end{equation}
Here, we have used (\ref{ohs_preconditioner}) and the fact that the decomposition $\bm{u}= P_0\bm{u} + (I-P_0) \bm{u}$ is $\bar{a}(\cdot, \cdot)$- orthogonal and $P_0$ is a projection; cf. \cite[Section 2.2]{toselli2005domain}. Because the range of $(I- P_0)$ is $\bm{V}_0^{\perp}$ , a subspace of $\bm{V}$, the Rayleigh quotient associated with the selfadjoint operator $(I-P_0)P_{OAS}(I-P_0)$ satisfies
$$
\lambda_{min}(P_{OAS}) \le
\min_{\bm{u} \in \bm{V}_0^{\perp} /\ \{\bm{0}\}}
\frac{\bar{a}(P_{OAS}(I-P_0)\bm{u}, (I-P_0)\bm{u})} {\bar{a}((I-P_0)\bm{u}, (I-P_0)\bm{u})}
$$
and
$$
\lambda_{max}(P_{OAS}) \ge
\max_{\bm{u} \in \bm{V}_0^{\perp} /\ \{\bm{0}\}}
\frac{\bar{a}(P_{OAS}(I-P_0)\bm{u}, (I-P_0)\bm{u})} {\bar{a}((I-P_0)\bm{u}, (I-P_0)\bm{u})}.
$$
We note that the maximum and minimum above  are considered on a subspace of $\bm{V}_h$. Then, from the Rayleigh quotient (\ref{Rayleigh_OHS}) and noting the $\bar{a}(\cdot, \cdot)$- orthogonal decomposition, we have
$$
\lambda_{min}(P_{OHS}) =
\min_{\bm{u} \ne \{\bm{0}\}}
\frac{\bar{a}(P_0\bm{u}, P_0\bm{u})+\bar{a}(P_{OAS}(I-P_0)\bm{u}, (I-P_0)\bm{u})} {\bar{a}(P_0\bm{u}, P_0\bm{u})+\bar{a}((I-P_0)\bm{u}, (I-P_0)\bm{u})} \ge \lambda_{min}(P_{OAS})
$$
and
$$
\lambda_{max}(P_{OHS}) =
\max_{\bm{u} \ne \{\bm{0}\}}
\frac{\bar{a}(P_0\bm{u}, P_0\bm{u})+\bar{a}(P_{OAS}(I-P_0)\bm{u}, (I-P_0)\bm{u})} {\bar{a}(P_0\bm{u}, P_0\bm{u})+\bar{a}((I-P_0)\bm{u}, (I-P_0)\bm{u})} \le \lambda_{max}(P_{OAS}).
$$
Therefore, the desired conclusion holds.

By combining the conclusion of Theorem \ref{th:est_AS2level} and the conclusion in Lemma \ref{th:ohs_comp_oas}, we have the following  condition number estimate of the overlapping hybrid Schwarz operator.

\begin{thm}\label{th:est_HS2level}
Assuming that the same assumptions of Theorem \ref{th:est_AS2level} hold, we have
\begin{equation}\label{ohs-cond}
\kappa(P_{OHS}) \leq C(H/\delta)^3(1+\log(H/\delta))(1+\log(H/h)).
\end{equation}
\end{thm}

The convergence rate of the overlapping multiplicative Schwarz preconditioner for the symmetric positive definite reformulated operator can be characterized by estimating the error propagation operator $E_{OMS}$. Defining
$$
||E_{OMS}||_{\bar{a}}^2 :=
\sup_{\bm{v} \in \bm{V}^h}
\frac{\bar{a}(E_{OMS}\bm{v}, E_{OMS}\bm{v})}{\bar{a}(\bm{v},\bm{v})},
$$
we then have

\begin{thm} \label{oms_est}
The error propagation operator of the overlapping multiplicative Schwarz method satisfies
$$
||E_{OMS}||_{\bar{a}}^2= ||I - P_{OMS}||_{\bar{a}}^2 \le 1- \frac{1}{3C_0^2} <1,
$$
where the constant $C_0$ is defined in (\ref{C-zero}).
\end{thm}

For the proof, one can check the assumptions for the multiplicative Schwarz theory in Chapter 2 of \cite{toselli2005domain}. We only need to note that the stable decomposition proof has been given in our previous work \cite{cai2015overlapping} for the overlapping additive Schwarz operator; the constant in the local stability estimate is 1 as we use exact solvers for subdomain problems.

For more discussions on the convergence rate of multiplicative Schwarz method, we refer the readers to
\cite{griebel1995abstract}. Some recent progresses on this topic can be found in \cite{nabben2003comparisons, notay2013further}. In our numerical tests, we accelerate OMS with the GMRES method since the preconditioner is not symmetric, and we find as expected that the two-level multiplicative Schwarz preconditioner has a much faster convergence rate than the additive or hybrid preconditioners. Alternatively, the symmetrized version of the multiplicative Schwarz preconditioner could be accelerated by PCG.



\section{Overlapping Schwarz preconditioners for the saddle point formulation}

Different from the approach in the previous section, one can directly apply the GMRES method as the outer iteration for the indefinite system (\ref{mixed_el}). Our goal of this section is to introduce some indefinite Schwarz preconditioners directly constructed from the global saddle point systems;  see the related works \cite{klawonn1998overlapping, klawonn2000comparison, pavarino1998preconditioned}.


In our construction, the local displacement (or velocity) space is always set to be
$$
\bm{V}_{i} = \bm{V}^h \cap (H^1_0 (\Omega_i^\prime))^d.
$$
The local pressure spaces can be one of following three choices.
\begin{itemize}

\item{\it Version 1.} We first impose the traditional zero mean value constraint for the  pressure  in $\Omega_i^\prime$, so that the associated pressure subspaces are
$$
\ds
Q_{i} = \{ Q^h \in L^2_0(\Omega_i^\prime) \}, \quad i=1, ..., N.
$$

\item{\it Version 2.}  In addition to imposing the zero mean value constraint for the  pressure  in $\Omega_i^\prime$, we also
set to zero the pressure degrees of freedom in the elements $K$ that touch $\partial \Omega^{'}/\ \partial \Omega$. That is, the associated pressure spaces are
$$
\ds
Q_{i} = \{ Q^h \in L^2_0(\Omega_i^\prime) : Q^h|_K=0, ~\forall K: \overline{K}\cap(\partial \Omega_i^\prime/ \partial \Omega) \neq \emptyset \}, \quad i=1, ..., N.
$$

\item{\it Version 3.} The pressure subspaces are same as {\it Version 2} except that the local zero mean value constraint for the pressures is removed. Therefore, we have
$$
\ds
Q_{i} = \{ Q^h \in L^2(\Omega_i^\prime) : Q^h|_K=0, ~\forall K: \overline{K}\cap(\partial \Omega_i^\prime/ \partial \Omega) \neq \emptyset \}, \quad i=1, ..., N.
$$

\end{itemize}

We remark that the pressure subspaces in {\it Version 2} are smaller than those in {\it Version 1}. In {\it Version 3}, the construction of pressure subspaces is a combination of {\it Version 1} and {\it Version 2}, see also  \cite{klawonn1998overlapping, pavarino1998preconditioned, pavarino2000indefinite}. In {\it Version 3}, the discrete local problems are nonsingular because setting to zero the pressures in the boundary elements of each subdomain has the effect of enforcing a zero Dirichlet boundary condition for the local pressures, which guarantees a unique pressure solution. This version is cheaper than {\it Version 2} since there is no need to enforce the zero mean value constraint in each local solve, but on the other hand the local problems in the finite element case are no longer guaranteed to satisfy a uniform inf-sup condition.

A coarse saddle point space $\bm{V}_0\times Q_0$ is obtained by using on the coarse mesh $\tau_H$  ${\bf Q}_2^H$ discontinuous $P_1^H$ mixed finite elements or
 ${\bf Q}_2$ discontinuous $Q_0$ mixed spectral elements, according to the discretization considered.

Given local and coarse (for $i=0$) pressure embedding operators ${{R}^p_{i}}^T:Q_i\rightarrow Q^h $, $i=0, 1,..,N$, we
can then decompose the discrete space $\bm{V}^h\times Q^h$ into local and coarse spaces as
\[
\bm{V}^h\times Q^h = \sum_{i=0}^N({R}_{i}^T\bm{V}_{i} \times {{R}_{i}^p}^TQ_i).
\]

Define the local (for $i \geq 1$) operators 
$\bar{P}_{i}^{mix}=
\left[\begin{array}{c}
\bar{{P}}_{i}^u\\
\bar{{P}}_{i}^p\\
\end{array}\right]
:\bm{V}^h\times Q^h\rightarrow \bm{V}_{i}\times Q_i$
by
\begin{equation}
\left\{
\begin{aligned}
& a_i^\prime(\bar{{P}}_{i}^u \bm{u}, \bm{v}) + b_i^\prime(\bm{v},\bar{{P}}_{i}^p p) = a(\bm{u}, {R}_{i}^T\bm{v}) + b({R}_{i}^T\bm{v}, p) \qquad \forall \bm{v} \in
\bm{V}_{i}, \\
&b_i^\prime(\bar{{P}}_{i}^u \bm{u},q) - c_i^\prime(\bar{{P}}_{i}^p p,q)= b(\bm{u},{{R}^p_{i}}^Tq) - c(p,{{R}^p_{i}}^Tq) \qquad \forall q \in Q_i,
\end{aligned}
\right.
\end{equation}
where the bilinear forms on the left hand sides are defined by integrals over $\Omega_i^\prime$. The operator for the coarse space is defined similarly.

Defining
${P}_{i}^{mix}:=
\left[\begin{array}{c}
{P}_{i}^u\\
{P}_{i}^p\\
\end{array}\right]
=
\left[\begin{array}{c}
{R}_{i}^T\bar{P}_{i}^u\\
{{R}^p_{i}}^T\bar{P}_{i}^p\\
\end{array}\right]$,
$i=0, 1,..,N$, our two-level Overlapping Additive Schwarz (OAS) operator formally has the same structure as before
\begin{equation}\label{ADDSCH_MIXED}
P_{OAS}^{mix}={P}_0^{mix}+\sum_{i=1}^N {P}_{i}^{mix}.
\end{equation}
Its operator form is
$
P_{OAS}^{mix}={B_{OAS}^{mix}}{{A}}^{mix},
$
with the mixed operator ${{A}}^{mix}$ obtained by subassembling the local operators (\ref{local_saddle_pt}) and with the mixed preconditioner
\begin{equation}
{B}_{OAS}^{mix}
= \left[\!\!\begin{array}{c}
R_0^{T}
{{R}_0^p}^{T}\\
\end{array}\!\!\right]
{A^{mix}_0}^{-1}
\left[\!\!\begin{array}{c}
R_0\\
{{R}_0^p}\\
\end{array}\!\!\right]
 +
 \sum_{i=1}^N
 \left[\!\!\begin{array}{c}
R_{i}^{T}
{{R}^p_{i}}^{T}\\
\end{array}\!\!\right]
\left[\!\!\begin{array}{cc}
\mu_i A^\prime_{i} & {B^\prime_{i}}^T \\
B^\prime_{i} & -\frac{1}{\lambda_i} C^\prime_{i} \\
\end{array}\!\!\right]^{-1}
\left[\!\!\begin{array}{c}
R_{i}\\
{{R}^p_{i}}\\
\end{array}\!\!\right],
\label{B_OAS_mixed}
\end{equation}
where $A^{mix}_0=
\left[\begin{array}{c}
R_0\\
{{R}_0^p}\\
\end{array}\right]
A^{mix}
\left[\begin{array}{c}
R_0^{T}
{{R}_0^p}^{T}\\
\end{array}\right]$.
We remark that this preconditioner leads to a system with complex eigenvalues in spite of the symmetry of both the original system and the preconditioner. The symmetry cannot be recovered as in the case of Section 3 because now the preconditioner and the original system are both indefinite. Therefore, in general, we no longer can employ the conjugate gradient method but must resort to a more general Krylov space method such as GMRES.

As before, the hybrid Schwarz preconditioner for the saddle point problem is similar to the hybrid Schwarz algorithm in Section 3, and is obtained by simply replacing the single bilinear form $\bar{a}$ by the saddle point bilinear systems. The hybrid saddle point algorithm then gives the following preconditioned operator \cite{toselli2005domain}.
$$
P_{OHS}^{mix}=P_0^{mix} + (I-P_0^{mix})\sum_{i=1}^N P_i^{mix}(I-P_0^{mix}).
$$

The multiplicative Schwarz preconditioner is obtained analogously by solving in sequence all subdomain saddle point problems. The implementation is similar to the multiplicative Schwarz preconditioner discussed in Section 3, so we have the multiplicative preconditioned operator
$$
P_{OMS}^{mix}=I - (I-P_N^{mix})(I-P_{N-1}^{mix})\cdots(I-P_0^{mix}).
$$

Similar to the overlapping additive Schwarz method, we can write out the hybrid Schwarz preconditioner ${B}^{mix}_{OHS}$ and the multiplicative Schwarz preconditioner ${B}^{mix}_{OMS}$. For the saddle point formulation, when there is no coarse grid preconditioner, the additive preconditioner (\ref{ADDSCH_MIXED}) and the multiplicative Schwarz preconditioner degenerate to the corresponding one-level overlapping Schwarz preconditioners.

\section{Numerical experiments in the plane}

In this section, we report on results of numerical tests in 2D with the overlapping Schwarz preconditioners for both the almost incompressible elasticity system and the incompressible Stokes problem. Although the focus of this work is hybrid and multiplicative Schwarz algorithms, the results based on the additive Schwarz algorithm will also be reported when necessary. Our problem is discretized with ${\bf Q}_2^h-P_1^h$ mixed finite elements or spectral elements. In Section 5.1, we report the numerical results for the positive definite reformulation of the mixed linear elasticity problem, while in Section 5.2, we report the results for the saddle point formulation. In Section 5.3, we show the robustness of the two-level overlapping Schwarz algorithms when the Poisson ratio approaches the incompressible limit and when jumps in the Lam\'{e} parameters are present. In Section 5.4, numerical results based on spectral elements are reported.

The domain is decomposed into $N$ overlapping subdomains of characteristic size $H$ and the overlap size $\delta$ is the minimal thickness of the extension $\Omega_i^\prime \setminus \Omega_i$ of each subdomain $\Omega_i$. If the resulting linear system is symmetric positive definite, we use PCG method, otherwise we use the GMRES. In all cases, we use one- or two-level overlapping Schwarz preconditioners as defined in Section 3 and Section 4. The initial guess is zero and the stopping criterion is set as a $10^{-6}$ reduction of the residual norm. In each test, we report the iteration counts (it.), the iteration errors (err.), i.e., the difference between the iterative solution and the solution obtained by using a direct solver. In the positive definite case, we also report the condition number (cond.) of the preconditioned operator defined as the ratio of its extreme eigenvalues $\lambda_{max}/\lambda_{min}$.

\subsection{Symmetric positive definite reformulation of mixed linear elasticity}


\begin{table}[htbp]
\begin{center}
  \centering
  \begin{tabular}{c|ccc|ccc}
    \hline
    \multicolumn{7}{c}{pos. def. reformulation, PCG-OHS(2)} \\
&\multicolumn{3}{c|}{$\delta = h$} &\multicolumn{3}{c}{$\delta = 2h$}\\
         {$N$}         &it.  &err.  & cond.=$\lambda_{max}/\lambda_{min}$ &it.  &err.  & cond.=$\lambda_{max}/\lambda_{min}$ \\
    \hline
    $2\times2$ &36 &1.1e-7 &153.5=4.000/2.61e-2  &24 &8.4e-8 &39.96=3.999/1.00e-1  \\
    $3\times3$ &59 &1.0e-6 &131.9=4.000/3.03e-2  &35 &1.0e-7 &37.37=4.000/1.07e-1  \\
    $4\times4$ &75 &4.5e-7 &121.0=4.000/3.31e-2  &40 &5.3e-7 &33.16=4.000/1.21e-1  \\
    $5\times5$ &75 &7.1e-7 &119.6=4.000/3.34e-2  &42 &5.2e-7 &33.84=4.000/1.15e-1  \\
    $6\times6$ &74 &1.1e-6 &116.7=4.000/3.43e-2  &41 &7.1e-7 &32.63=4.000/1.23e-1  \\
    $7\times7$ &74 &1.6e-6 &117.4=4.000/3.41e-2  &41 &7.3e-7 &32.17=4.000/1.24e-1  \\
    $8\times8$ &76 &7.9e-7 &111.9=4.000/3.58e-2  &41 &1.1e-6 &32.48=4.000/1.23e-1  \\
    \hline
  \end{tabular}
\caption{Scalability of PCG-OHS(2) for the positive definite reformulation. Iteration counts, errors, condition numbers and extreme eigenvalues for increasing number of subdomains $N$. Fixed $H/h=9,\ \nu=0.4999$.
}\label{scalability-pcg-ohs}
\end{center}
\end{table}

\begin{figure}[h]
\centering
\vspace{-0.8cm}
\includegraphics[scale=0.45]{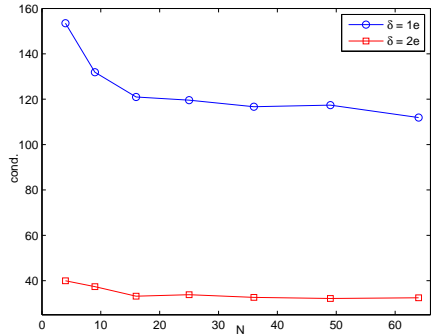}
\caption{Plot of cond. from Table \ref{scalability-pcg-ohs}.}
\label{cond_num_ohs}
\end{figure}

\begin{table}[h]
\begin{center}
  \centering
  \begin{tabular}{c|cc|cc|cc|cc|cc|cc}
    \hline
     &\multicolumn{4}{c|}{pos. def, GMRES-OMS(1)} &\multicolumn{4}{c|}{pos. def, GMRES-OMS(2)} &\multicolumn{4}{c}{pos. def, GMRES-OAS(2)}  \\
  &\multicolumn{2}{c|}{$\delta = h$} &\multicolumn{2}{c|}{$\delta = 2h$} &\multicolumn{2}{c|}{$\delta = h$} &\multicolumn{2}{c|}{$\delta = 2h$} &\multicolumn{2}{c|}{$\delta = h$} &\multicolumn{2}{c}{$\delta = 2h$}\\
       {$N$}    &it. &err. &it.  &err. &it. &err. &it.  &err.     &it.  &err.  &it.  &err. \\
    \hline
    $2\times2$ &10 &1.2e-7  &6  &2.0e-7 &9  &1.9e-8  &5  &3.3e-7  &24 &2.6e-6  &16 &7.6e-7  \\
    $3\times3$ &17 &2.2e-6  &10 &2.5e-7 &14 &1.1e-6  &8  &8.8e-7  &39 &4.4e-6  &22 &7.8e-7  \\
    $4\times4$ &29 &9.4e-6  &14 &1.5e-6 &16 &2.8e-6  &9  &2.1e-6  &44 &8.7e-6  &26 &2.3e-6  \\
    $5\times5$ &40 &2.1e-5  &18 &9.2e-7 &18 &2.6e-6  &10 &1.0e-6  &43 &9.4e-6  &25 &2.8e-6  \\
    $6\times6$ &58 &2.0e-5  &26 &4.0e-6 &17 &1.1e-5  &10 &1.6e-6  &45 &2.1e-6  &24 &4.7e-6  \\
    \hline
  \end{tabular}
\caption{Scalability of OMS(2) for the symmetric positive definite reformulation. Iteration counts and errors for increasing number of subdomains $N$. Fixed $H/h=5,\ \nu=0.4999$.
}\label{scalability-gmres-oms}
\end{center}
\end{table}

{\bf Scalability of OHS(2) and OMS(2).}
We investigate the scalability of the hybrid and multiplicative overlapping Schwarz preconditioners proposed in Section 3.

In Table \ref{scalability-pcg-ohs}, we consider the two-level PCG-OHS(2) algorithms, applied to the positive definite reformulation, for increasing number of subdomains $N$ and fixed $H/h=9, \nu=0.4999$. 
In each case, we consider both a minimal overlap of $\delta = h$ and a larger overlap of $\delta = 2h$. The results show that the PCG-OHS(2) iteration count is bounded from above by a constant independent of $N$, clearly showing the scalability of the proposed preconditioners. Increasing the overlap size yields a considerable improvement for the positive definite reformulation. The condition numbers from this table are also plotted in Figure \ref{cond_num_ohs} as a function of $N$. All the results clearly show the scalability of PCG-OHS(2). 

\begin{table}
\begin{center}
  \centering
  \begin{tabular}{r|cr|rr}
    \hline
\multicolumn{5}{c}{pos. def. reformulation, PCG-OHS(2)} \\
 &\multicolumn{2}{c|}{$\nu = 0.3$} &\multicolumn{2}{c}{$\nu = 0.4999$}\\
          {$H/\delta$}       &it.   & cond.=$\lambda_{max}/\lambda_{min}$  &it.    & cond.=$\lambda_{max}/\lambda_{min}$  \\
    \hline
 $5.33$   &14   &4.75=4.000/8.42e-1      &25   &57.97=4.000/6.90e-2  \\  
 $6.40$   &15   &5.25=4.000/7.62e-1      &30   &84.63=4.000/4.73e-2   \\  
 $8.00 $  &16   &6.63=4.000/6.03e-1      &35   &127.6=4.000/3.14e-2  \\  
 $10.67$  &18   &8.58=4.000/4.66e-1      &42   &191.9=4.000/2.08e-2  \\  
  $16.00$ &20  &12.03=4.000/3.33e-1      &59   &268.3=4.000/1.49e-2  \\  
 $21.33$  &22  &15.45=4.000/2.59e-1      &67   &293.6=4.000/1.36e-2  \\  
 $32.00$  &25  &22.26=4.000/1.80e-1     &112   &840.3=4.000/4.76e-3  \\  
 $42.67$  &29  &29.07=4.000/1.38e-1     &146  &1814.3=4.000/2.21e-3  \\  
 $64.00$  &33  &42.63=4.000/9.38e-2     &224  &4964.0=4.000/8.06e-4 \\ 
 $128.00$ &45  &83.32=4.000/4.80e-2     &466  &22883.5=4.000/1.75e-4  \\ 
    \hline
  \end{tabular}
  \caption{$H/\delta$-dependence of PCG-OHS(2), positive definite reformulation. Iteration counts, condition numbers and extreme eigenvalues for increasing $H/\delta$. Fixed $N= 2\times 2, H/h=128$,  $\nu = 0.3$ (left), $\nu=0.4999$ (right).
}\label{pcg-ohs-Hdelta_04999}
\end{center}
\end{table}

An analogous scalability test is reported in Table \ref{scalability-gmres-oms} for GMRES-OMS(2) applied to the positive definite reformulation with fixed $H/h=5, \nu=0.4999$, overlap $\delta = h$ (left) and $\delta = 2h$ (right), for $N$ increasing up to $6\times 6$ subdomains. The results show that the GMRES-OMS(2) iteration counts are bounded from above by a constant independent of $N$, clearly showing the scalability of the two-level overlapping multiplicative Schwarz preconditioner, while the one-level algorithm GMRES-OMS(1) is not scalable since its iteration counts grow with $N$. We also report the numerical results based on the GMRES-OAS(2) in  Table \ref{scalability-gmres-oms}. By comparing the number of iterations needed for GMRES-OMS(2) and those based on GMRES-OAS(2), we see that GMRES-OAS(2) needs more than two times of number of iteration counts of GMRES-OMS(2), i.e.  GMRES-OMS(2) converges about twice as fast as GMRES-OAS(2).

\begin{figure}[t]
\centerline{
\vspace{-0.8cm}
\includegraphics[height=5.5cm]{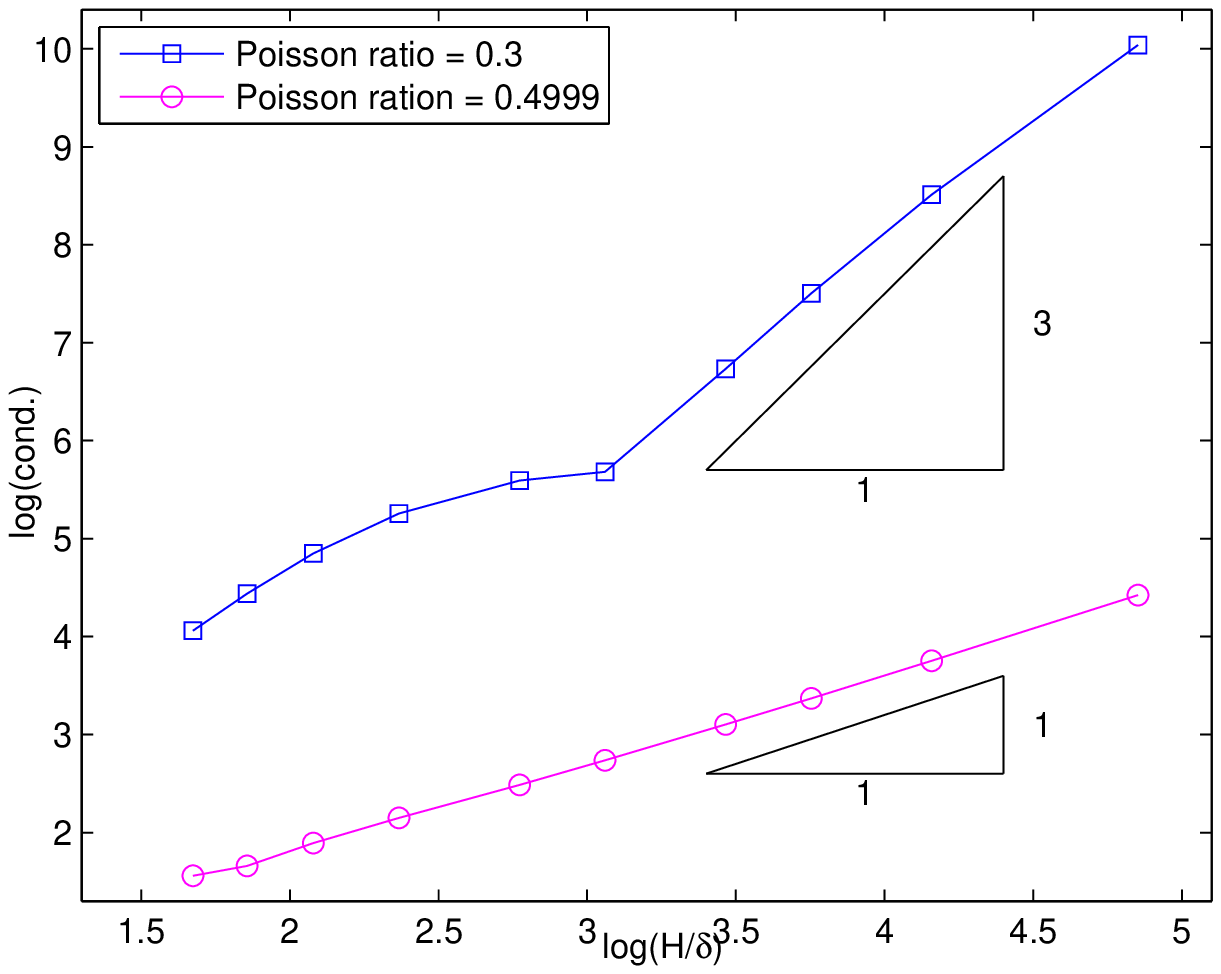}
}
\caption{Plot of PCG-OHS(2) condition numbers from Table \ref{pcg-ohs-Hdelta_04999}.
}
\label{compare_OHS_diff_Hdelta}
\end{figure}

{\bf OHS(2) dependence on $H/\delta$.}
In order to check our main bound in Theorem \ref{th:est_HS2level} that predicts a $(H/\delta)^3$ growth of the condition number for the almost incompressible case, we have investigated the effect of increasing the ratio $H/\delta$ while fixing $N = 4$ and $H/h = 128$. The results reported in Table \ref{pcg-ohs-Hdelta_04999}, and also plotted in Figure \ref{compare_OHS_diff_Hdelta}, confirm the theoretical $(H/\delta)^3$ bound in the almost incompressible case, while the bound appears to be only linear in $H/\delta$ in the compressible case.

{\bf OHS(2) and OAS(2) dependence on $H/h$.}
We then investigate the OAS(2) and OHS(2) dependence on the ratio $H/h$ for the positive definite reformulation, considering both the compressible and almost incompressible cases. Numerical results are reported in Table \ref{pcg-oas-ohs-Hoverh} for both the OAS(2) and the OHS(2) tests. With fixed $N = 2\times 2$ and relative overlap size $H/\delta = 2$, we increase the ratio $H/h$ from 4 to 64. The results indicate that the condition numbers for both the compressible and almost incompressible cases seem to be independent of $H/h$. Moreover, we observe that both the maximum and minimum eigenvalues of PCG-OAS(2) and PCG-OHS(2) seem to be independent of $H/h$, if the value of $H/\delta$ is fixed. In comparison, the condition number of the overlapping Schwarz algorithms developed in \cite{dohrmann2009overlapping, dohrmann2010hybrid} increases linearly as $1+\log(H/h)$ increases and $H/\delta$ is fixed (see Table 7.2 and Figure 7.2 in \cite{dohrmann2009overlapping}, Table 2 and Figure 2 in \cite{dohrmann2010hybrid}, for both 2D and 3D configurations). We do not have a very precise explanation for this difference, which could be due to the use of different coarse problems and it might not hold for the more general 3D problems considered in \cite{dohrmann2009overlapping}. Nevertheless, in our particular configuration, it seems that it could be possible to improve the bounds (\ref{oas-cond}) and (\ref{ohs-cond}) to
$$
\kappa(P_{OAS}) \leq C(H/\delta)^3(1+\log(H/\delta)) \quad \mbox{and} \quad \kappa(P_{OHS}) \leq C(H/\delta)^3(1+\log(H/\delta)),
$$
indicating that the two-level overlapping Schwarz algorithms with standard coarse spaces developed in this work are actually optimal, i.e., in the generous overlap case $H/\delta=Constant$, these bounds become independent of $H/h$.
Moreover, we also tested our preconditioners for the checkerboard test (see the coefficient distribution in Section 5.3) in which neither Assumption 1 nor Assumption 2 are satisfied. Again, the resulting condition numbers (not reported here) seem to be constants if $H/\delta$ is fixed, indicating that the {\it quasi-monotonicity} condition can be weaken or removed, at least  under the special configurations considered in this work.

\begin{table}[h]
\begin{center}
  \centering
  \begin{tabular}{r|ccl|ccl||ccl|ccl}
  \hline
  &\multicolumn{6}{c||}{pos. def. reformulation, PCG-OAS(2)} & \multicolumn{6}{c}{pos. def. reformulation, PCG-OHS(2)}\\
  &\multicolumn{3}{c|}{$\nu = 0.3$} &\multicolumn{3}{c||}{$\nu = 0.4999$} &\multicolumn{3}{c|}{$\nu = 0.3$} &\multicolumn{3}{c }{$\nu = 0.4999$} \\
{$H/h$} &it. &err. &cond. &it.  &err.  &cond. &it. &err. &cond. &it. &err. &cond.\\
  \cline{1-13}
   $4$ &14 &8.5e-7 &5.19  &24 &3.8e-8 &38.39 &12 &8.6e-7 &4.33  &22 &5.9e-8 &30.69 \\
   $8$ &14 &2.5e-6 &5.16  &24 &6.8e-8 &38.42 &13 &1.2e-6 &4.37  &22 &6.2e-8 &30.73  \\
  $16$ &15 &4.3e-6 &5.16  &24 &1.9e-7 &38.42 &12 &6.0e-6 &4.30  &22 &2.0e-7 &30.73 \\
  $32$ &15 &1.5e-5 &5.17  &23 &1.1e-6 &38.43 &12 &2.4e-5 &4.31  &21 &1.5e-6 &30.73 \\
  $64$ &15 &3.2e-5 &5.17  &23 &2.4e-6 &38.43 &13 &4.1e-5 &4.30  &22 &1.3e-6 &30.73  \\
    \hline
  \end{tabular}
\caption{$H/h$-dependence of PCG-OAS(2) and PCG-OHS(2), positive definite reformulation. Iteration counts, errors, condition numbers for increasing $H/h$. Fixed $H/\delta= 4,\ N=2 \times 2$.
}\label{pcg-oas-ohs-Hoverh}
\end{center}
\end{table}

\subsection{Saddle point preconditioners for the Stokes system}

\begin{table}[h]
\begin{center}
  \centering
  \begin{tabular}{c|cc|cc|cc||cc|cc|cc}
    \hline
    &\multicolumn{6}{c||}{saddle point formulation, GMRES-OAS(1)}  &\multicolumn{6}{c}{saddle point formulation, GMRES-OAS(2)}  \\
      &\multicolumn{2}{c|}{V1} &\multicolumn{2}{c|}{V2} &\multicolumn{2}{c||}{V3} &\multicolumn{2}{c|}{V1} &\multicolumn{2}{c|}{V2} &\multicolumn{2}{c}{V3} \\
    \cline{1-13}
     N    &it. &err.  &it. &err.  &it. &err. &it. &err. &it. &err. &it. &err. \\
    \hline
    $2\times2$ &21  &2.6e-5  &15 &6.2e-6 &14  &9.8e-6 &20 &3.6e-5 &16 &7.4e-6 &16  &1.3e-5 \\
    $3\times3$ &35  &1.7e-5  &20 &2.2e-5 &20  &2.0e-5 &23 &3.2e-5 &17 &2.1e-5 &18  &2.6e-5 \\
    $4\times4$ &59  &2.6e-4  &24 &3.7e-5 &32  &6.2e-5 &25 &4.4e-5 &18 &1.4e-5 &19  &3.6e-5 \\
    $5\times5$ &83  &1.8e-3  &29 &4.8e-5 &42  &1.0e-4 &25 &1.4e-4 &18 &3.2e-5 &20  &3.7e-5 \\
    $6\times6$ &112 &5.6e-3  &31 &7.6e-5 &51  &1.9e-4 &25 &2.3e-4 &18 &2.8e-5 &20  &4.7e-5 \\
    \hline
  \end{tabular}
\vspace{2mm} \caption{GMRES-OAS for Stokes problem with different versions of pressure local spaces. Iteration counts and errors for increasing number of subdomains $N$. Fixed $H/h=5$, $\delta=h$.
}\label{Stokes_OAS_1e}
\end{center}
\end{table}

\begin{table}[h]
\begin{center}
  \centering
  \begin{tabular}{c|cc|cc|cc||cc|cc}
    \hline
    &\multicolumn{6}{c||}{saddle point formulation, GMRES-OMS(1)}  &\multicolumn{4}{c}{saddle point formulation}  \\
      &\multicolumn{2}{c|}{V1} &\multicolumn{2}{c|}{V2} &\multicolumn{2}{c||}{V3} &\multicolumn{2}{c|}{GMRES-OMS(2), V2} &\multicolumn{2}{c}{GMRES-OHS(2), V2} \\
    \cline{1-11}
     N    &it. &err.  &it. &err.  &it. &err. &it. &err. &it. &err.  \\
    \hline
    $2\times2$ &11  &7.9e-6  &8  &4.9e-7 &7   &1.6e-6 &6  &1.1e-6 &15 &1.1e-5  \\
    $3\times3$ &19  &1.7e-5  &11 &4.6e-6 &12  &5.1e-6 &7  &6.6e-7 &16 &2.2e-5  \\
    $4\times4$ &28  &1.4e-4  &13 &2.4e-5 &17  &1.6e-5 &7  &2.0e-6 &16 &3.5e-5  \\
    $5\times5$ &39  &7.9e-5  &15 &3.7e-5 &24  &1.1e-5 &7  &7.1e-6 &16 &7.3e-5  \\
    $6\times6$ &51  &8.9e-4  &18 &4.7e-5 &29  &4.1e-5 &7  &9.8e-6 &16 &1.1e-4  \\
    \hline
  \end{tabular}
\vspace{2mm} \caption{Comparisons of GMRES-OMS(1), GMRES-OMS(2) and GMRES-OHS(2) for Stokes problem with different versions of pressure local spaces. Iteration counts and errors for increasing number of subdomains $N$. Fixed $H/h=5$, $\delta=h$.
}\label{Comp_OMS_OHS}
\end{center}
\end{table}

We now consider the Stokes problem with homogeneous Dirichlet boundary condition and compare different versions of the indefinite overlapping Schwarz preconditioners proposed in Sec. 4.

In Table \ref{Stokes_OAS_1e}, we report numerical results for both one-level and two-level overlapping additive Schwarz preconditioners. First, the number of iterations of the one-level algorithms  increases if the number of subdomains increases. Instead, the number of iterations of the two-level algorithms is bounded from above when the number of subdomains increases,  clearly showing the scalability of the two-level algorithms. Moreover, by comparing different versions of the algorithms, we see clearly that the algorithm based on {\it Version 2} gives the best performance.

In Table \ref{Comp_OMS_OHS}, we compare the performance of the indefinite multiplicative Schwarz preconditioners using different versions of pressure local spaces. In addition, the performance of GMRES-OHS(2) and GMRES-OMS(2) using the {\it Version 2} algorithm are reported. The results confirm that the algorithm based on {\it Version 2} is the best. Moreover, GMRES-OMS(2) gives much better convergence rate than GMRES-OHS(2).

\subsection{Robustness with respect to incompressibility and coefficient discontinuities}

In this section, we investigate the robustness of our overlapping Schwarz algorithms with respect to the discontinuities of
parameters and the Poisson ratios.

\begin{table}[h]
\begin{center}
  \centering

    \begin{tabular}{l|cc|cc|ccc|cc}
    \hline
   &\multicolumn{4}{c|}{saddle point formulation} &\multicolumn{5}{c}{pos. def. reformulation}  \\
   &\multicolumn{2}{c|}{GMRES-OHS(2)} &\multicolumn{2}{c|}{GMRES-OMS(2)} &\multicolumn{3}{c}{PCG-OHS(2)} &\multicolumn{2}{c}{GMRES-OMS(2)}  \\
     $\nu$  &it.   &err.    &it.   &err. &it.   &err.   & cond.=$\lambda_{max}/\lambda_{min}$  &it.   &err.      \\
    \hline
    0.4      &14  &2.3e-6  &5  &5.5e-6  &14 &4.1e-7  &\hphantom{1}4.69=3.996/8.53e-1   &4  &1.6e-5    \\
    0.49     &15  &3.6e-6  &6  &2.9e-6  &17 &3.2e-7  &\hphantom{1}6.22=3.999/6.43e-1   &7  &5.0e-7    \\
    0.499    &15  &4.4e-6  &6  &4.9e-6  &25 &1.2e-7  &15.79=4.000/2.53e-1              &10 &9.8e-7    \\
    0.4999   &15  &4.5e-6  &6  &5.1e-6  &31 &6.1e-8  &29.88=4.000/1.34e-1              &13 &5.1e-7    \\
    0.49999  &15  &4.5e-6  &6  &5.2e-6  &32 &5.2e-8  &38.44=4.000/1.04e-1              &13 &7.4e-7    \\
   0.499999  &15  &4.5e-6  &6  &5.2e-6  &33 &1.7e-8  &39.61=4.000/1.01e-1              &13 &7.7e-7    \\
    \hline
   0.5       &15  &2.9e-5  &6  &2.3e-6  & N/A  &    &   & N/A       &              \\
  \hline

  \end{tabular}
\caption{Robustness with respect to Poisson ratio for  PCG - OHS(2) and GMRES - OMS(2). Iteration counts, errors, condition numbers and extreme eigenvalues for an increasing Poisson ratio $\nu \rightarrow \frac{1}{2}$. Fixed $N = 3 \times 3,\ H/h = 4$, overlap $\delta = h$. 
}
\label{parameter_test}
\end{center}
\end{table}

{\bf Robustness of OHS(2) and OMS(2) for almost incompressible materials.}
In order to test the robustness of our algorithms for almost incompressible materials, we study the performance of different preconditioners under different Poisson ratios. The system is discretized with a fixed number $N = 3 \times 3$ of subdomains, with ratio $H/h = 4$, and overlap $\delta = h$. For the saddle point formulation, the system is solved by GMRES with either hybrid (OHS(2)) or multiplicative (OMS(2))  preconditioners. For the positive definite reformulation, we consider GMRES-OMS(2) or PCG-OHS(2). Table \ref{parameter_test} reports the iteration counts and the corresponding errors between iterative and direct solution when the Poisson ratio $\nu$ approaches 1/2. The results clearly show the robustness of both the OMS(2) and the OHS(2) preconditioners, with a better performance of the multiplicative OMS(2) preconditioner.
\begin{figure}
\centerline{
\includegraphics[height=4.6cm]{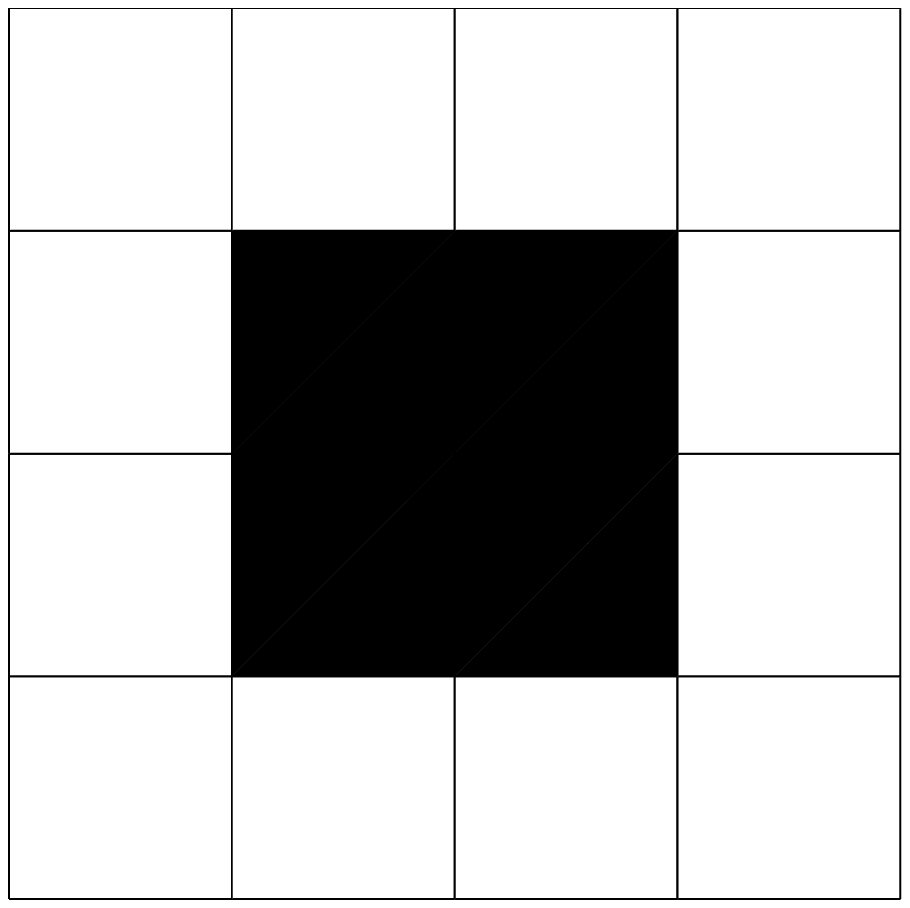}
\includegraphics[height=4.6cm]{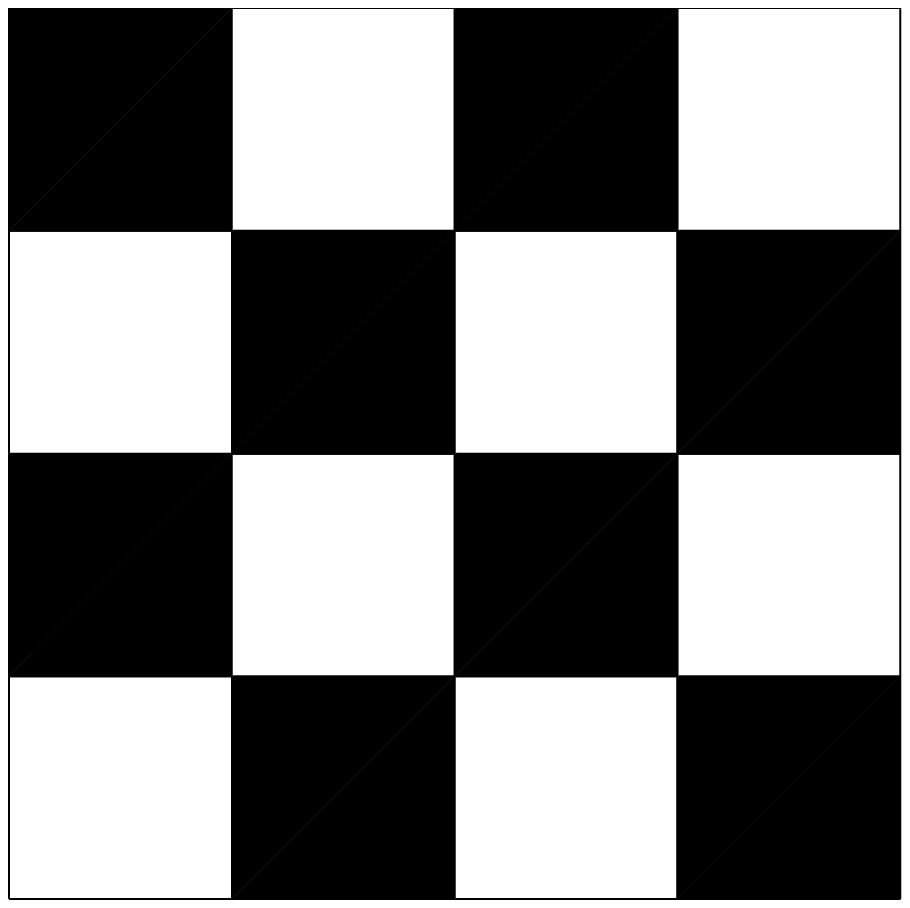}
}
\caption{
Central jump and checkerboard coefficient distributions.
}
\label{central_jump_checkerboard}
\end{figure}

\begin{table}[h]
\begin{center}
  \centering
  \begin{tabular}{l|cc|cc|ccc|cc}
    \hline
   &\multicolumn{4}{c|}{saddle point formulation} &\multicolumn{5}{c}{pos. def. reformulation}  \\
   &\multicolumn{2}{c|}{GMRES-OHS(2)} &\multicolumn{2}{c|}{GMRES-OMS(2)} &\multicolumn{3}{c|}{PCG-OHS(2)} &\multicolumn{2}{c}{GMRES-OMS(2)}  \\
      $\nu$  &it.   &err.    &it.   &err.   &it.   &err.  & cond.=$\lambda_{max}/\lambda_{min}$  &it.   &err.   \\
    \hline
       \cline{1-10}
    &\multicolumn{9}{c}{central jump test} \\
    0.3      &13 &2.5e-6  &5  &1.3e-6   &13 &1.1e-10  &4.44=3.997/9.00e-1  &4  &2.7e-9   \\
    0.4      &14 &1.6e-6  &5  &3.1e-6   &13 &1.3e-10  &4.48=3.997/8.92e-1  &4  &2.3e-9    \\
    0.49     &15 &4.8e-6  &5  &1.8e-5   &15 &6.7e-11  &5.30=3.995/7.54e-1  &5  &2.7e-9    \\
    0.499    &15 &6.8e-6  &5  &2.6e-5   &18 &6.0e-11  &7.21=3.997/5.54e-1  &6  &2.3e-9    \\
    0.4999   &15 &7.1e-6  &5  &2.8e-5   &20 &1.8e-11  &7.76=3.997/5.15e-1  &6  &3.9e-9     \\
    0.49999  &15 &7.2e-6  &5  &2.8e-5   &20 &2.3e-11  &7.83=3.997/5.11e-1  &6  &4.1e-9     \\
    \hline
        \cline{1-10}
    &\multicolumn{9}{c}{checkerboard test} \\
             &15 &5.2e-6  &6  &4.1e-6  &25 &2.0e-11   &8.86=3.998/4.51e-1  &8 &3.5e-9  \\
  \hline
  \end{tabular}
\caption{Robustness with respect to jumps of the elliptic coefficients. Saddle point and positive definite formulations using both OHS(2) and OMS(2) preconditioners. Iteration counts, errors, condition numbers and extreme eigenvalues for  increasing Poisson ratio $\nu \rightarrow \frac{1}{2}$. Fixed $N = 4 \times 4,\ H/h = 4$, overlap $\delta = h$. 
}
\label{jump_test}
\end{center}
\end{table}

{\bf Robustness of OHS(2) and OMS(2) with respect to discontinuous material parameters.}
We then consider two tests with discontinuous material parameters. The first, called "central jump" (left part of Figure \ref{central_jump_checkerboard}), consists of a square domain with $4\times 4$ subdomains, where the Poisson ratio $\nu$ equals the value given in the left column of the table in the $2\times 2$ central (black) subdomains, while $\nu = 0.3$ in the remaining (white) subdomains. In the second test, called "checkerboard" (right part of Figure \ref{central_jump_checkerboard}), the Poisson ratio is a piecewise constant function on each subdomain, with values varying randomly between $0.3$ and $0.49999$. In the checkerboard test, black subdomains are associated with almost incompressible materials, while blank subdomains with compressible material. More clearly, we set $E=6000$, while the Poisson ratio values $\nu$,  in the 4-by-4 subdomains, are given by
$$
\left[\begin{array}{cccc}
  0.49999 &0.37   &0.499  &0.41 \\
  0.3    &0.49999 &0.33   &0.4999\\
  0.49999 &0.29   &0.499  &0.3   \\
  0.2     &0.4999 &0.31   &0.499
\end{array}
\right].
$$
The coefficients $\mu$ and $\lambda$ are then calculated using (\ref{Poisson}). We remark that this coefficient distribution of $\mu$ does not satisfy Assumption 1 or Assumption 2.

In Table \ref{jump_test}, from left to right, the table reports the iteration counts and iteration errors of the GMRES-OHS(2) and GMRES-OMS(2) for the saddle point formulation, GMRES-OMS(2) and PCG-OHS(2) for the positive definite formulation, where for the latter case we report also the condition number and extreme eigenvalues. The results clearly show that for both the positive definite reformulation and the saddle point formulation, the proposed two-level multiplicative  and hybrid Schwarz algorithms are robust with respect to the jumps in the Poisson ratio in both the central jump and the checkerboard tests. Since the coefficient distribution of the checkerboard test actually does not satisfy Assumption 1 or Assumption 2, these results show that our algorithms still perform very well even when the quasi-monotonicity assumptions are not satisfied. It will be interesting to develop in future work a theoretical investigation of this numerical result.

\subsection{Results with $Q_n - Q_{n-2}$ spectral elements}

We finally tested our hybrid and multiplicative Schwarz preconditioners also on  $Q_n - Q_{n-2}$ spectral elements discretizations of almost incompressible elasticity and Stokes systems.
Table \ref{scalability-pcg-ohs_sem} shows the scalability in $N$ of PCG-OHS(2) for cubic ($n=3$) spectral elements for both minimal overlap $\delta = h$ and larger overlap $\delta = 2h$, while Table \ref{pcg-oas-ohs-Hoverh_sem} shows the optimality in $H/h$ of both the additive and hybrid preconditioners. The results of Table  \ref{pcg-oas-ohs-Hoverh_sem} are strikingly close to the quadratic finite element results of Table \ref{pcg-oas-ohs-Hoverh}.

Tables \ref{parameter_test_sem} and \ref{jump_test_sem} show the robustness of OHS(2) and OMS(2) with respect to the Poisson ratio $\nu$ approaching 1/2 and to discontinuous material parameters across subdomain interfaces, for both the saddle point and positive definite formulations.

The last Table \ref{degree-pcg-ohs-oms_sem} shows the $n$-independence of PCG-OHS(2) and GMRES-OMS(2) for the positive definite reformulation. All the iteration counts and condition numbers are basically independent of the $Q_n - Q_{n-2}$ polynomial degree when it increases from 3 to 8. This is due to the fact that the overlap size $\delta = h$ corresponds to a whole spectral element, hence the overlap is generous in terms of GLL points, that play the role of the fine mesh in the spectral element case.
In summary, these results confirm that the good convergence properties of our hybrid and multiplicative Schwarz preconditioners hold for both finite element and spectral element discretizations.

\begin{table}[t]
\begin{center}
 \centering
 \begin{tabular}{c|ccc|ccc}
 \hline
 \multicolumn{7}{c}{pos. def. reformulation, PCG-OHS(2)} \\
&\multicolumn{3}{c|}{$\delta = h$} &\multicolumn{3}{c}{$\delta = 2h$}\\
 {$N$} &it. &err. & cond.=$\lambda_{max}/\lambda_{min}$ &it. &err. & cond.=$\lambda_{max}/\lambda_{min}$ \\
 \hline
 $2\times2$ &26 &4.0e-8 &50.28=4.000/7.95e-2 &18 &1.1e-7 &11.04=4.000/3.62e-1 \\
 $3\times3$ &38 &1.8e-8 &45.00=4.000/8.89e-2 &22 &2.1e-7 &11.84=4.000/3.38e-1 \\
 $4\times4$ &47 &2.1e-7 &40.90=4.000/9.78e-2 &25 &3.2e-7 &11.49=4.000/3.48e-1 \\
 $5\times5$ &48 &4.4e-7 &42.17=4.000/9.48e-2 &25 &4.1e-7 &10.62=4.000/3.76e-1 \\
 $6\times6$ &48 &4.0e-7 &40.09=4.000/9.98e-2 &25 &2.5e-7 &10.77=4.000/3.71e-1 \\
 $7\times7$ &48 &7.9e-7 &38.40=4.000/1.01e-1 &26 &2.9e-7 &10.80=4.000/3.70e-1 \\
 $8\times8$ &46 &8.7e-7 &38.54=4.000/1.04e-1 &25 &4.4e-7 &10.45=4.000/3.83e-1 \\
 \hline
 \end{tabular}
\caption{$Q_n-Q_{n-2}$ spectral elements. Scalability of PCG-OHS(2) for the positive definite reformulation. Iteration counts, errors, condition numbers and extreme eigenvalues for increasing number of subdomains $N$. Fixed $n=3$,
$H/h=5,\ \nu=0.4999$. 
}\label{scalability-pcg-ohs_sem}
\end{center}
\end{table}

\begin{table}[h]
\begin{center}
  \centering
  \begin{tabular}{r|ccl|ccl||ccl|ccl}
  \hline
  &\multicolumn{6}{c||}{pos. def. reformulation, PCG-OAS(2)} & \multicolumn{6}{c}{pos. def. reformulation, PCG-OHS(2)}\\
  &\multicolumn{3}{c|}{$\nu = 0.3$} &\multicolumn{3}{c||}{$\nu = 0.4999$} &\multicolumn{3}{c|}{$\nu = 0.3$} &\multicolumn{3}{c }{$\nu = 0.4999$} \\
{$H/h$} &it. &err. &cond. &it.  &err.  &cond. &it. &err. &cond. &it. &err. &cond.\\
  \cline{1-13}
   $4$  &14 &1.0e-6   &5.24  &22 &3.1e-7 &35.60 &13 &5.4e-7 &4.53  &20 &1.2e-7 &27.99 \\
   $8$  &14 &2.1e-6   &5.33  &24 &8.7e-8 &37.35 &13 &1.1e-6 &4.55  &21 &3.3e-7 &29.78  \\
  $16$  &14 &6.7e-6   &5.28  &22 &8.4e-8 &38.08 &13 &1.9e-6 &4.54  &20 &7.6e-7 &30.43 \\
  $32$  &13 &1.7e-5   &5.16  &24 &2.9e-7 &38.32 &12 &1.9e-5 &4.46  &22 &6.9e-7 &30.64 \\
  $64$  &13 &3.5e-5   &5.32  &24 &1.1e-6 &38.40 &12 &2.5e-5 &4.55  &22 &1.8e-6 &30.71  \\
    \hline
  \end{tabular}
\caption{$Q_n - Q_{n-2}$ spectral elements with $n=2$. $H/h$-dependence of PCG-OAS(2) and PCG-OHS(2), positive definite reformulation. Iteration counts, errors, condition numbers and extreme eigenvalues for increasing $H/h$. Fixed $H/\delta= 4,\ N=2 \times 2$.
}\label{pcg-oas-ohs-Hoverh_sem}
\end{center}
\end{table}

\begin{table}[h]
\begin{center}
  \centering
  \begin{tabular}{l|cc|cc|ccc|cc}
    \hline
   &\multicolumn{4}{c|}{saddle point formulation} &\multicolumn{5}{c}{pos. def. reformulation}  \\
   &\multicolumn{2}{c|}{GMRES-OHS(2)} &\multicolumn{2}{c|}{GMRES-OMS(2)} &\multicolumn{3}{c|}{PCG-OHS(2)} &\multicolumn{2}{c}{GMRES-OMS(2)}  \\
     $\nu$  &it.   &err.    &it.   &err.   &it.   &err.   & cond.=$\lambda_{max}/\lambda_{min}$ &it.   &err.      \\
    \hline
  0.4       &15  &1.1e-6  &5  &1.1e-6  &14 &7.1e-7 &\hphantom{1}4.73=3.99/0.84  &5  & 1.1e-6  \\
  0.49      &17  &2.6e-6  &7  &1.0e-6  &16 &1.0e-6 &\hphantom{1}6.18=3.99/0.65  &7  & 5.8e-7  \\
  0.499     &18  &1.9e-6  &7  &3.4e-6  &24 &2.7e-7 &13.91=4.0/0.29              &10 & 4.3e-7  \\
  0.4999    &18  &2.1e-6  &8  &1.7e-7  &28 &1.1e-7 &26.27=4.00/0.15             &13 & 3.8e-7  \\
  0.49999   &18  &2.1e-6  &8  &1.7e-7  &29 &1.2e-7 &35.28=4.00/0.11             &13 &  5.9e-7 \\
  0.499999  &18  &2.1e-6  &8  &1.7e-7  &29 &1.2e-7 &36.60=4.00/0.109            &13 &  6.1e-7 \\
    \hline
  0.5       &18  &2.1e-6  &8  &1.7e-7  &N/A  &       &      &  N/A      &              \\
  \hline
  \end{tabular}
\caption{$Q_n - Q_{n-2}$ spectral elements with $n=2$.  Robustness with respect to the Poisson ratio for PCG - OHS(2) and GMRES - OMS(2). Iteration counts, errors, condition numbers and extreme eigenvalues for increasing Poisson ratio $\nu \rightarrow \frac{1}{2}$. Fixed $N = 3 \times 3,\ H/h = 4$, overlap $\delta = h$. 
}
\label{parameter_test_sem}
\end{center}
\end{table}

\begin{table}[h]
\begin{center}
  \centering
  \begin{tabular}{l|cc|cc|ccc|cc}
    \hline
   &\multicolumn{4}{c|}{saddle point formulation} &\multicolumn{5}{c}{pos. def. reformulation}  \\
   &\multicolumn{2}{c|}{GMRES-OHS(2)} &\multicolumn{2}{c|}{GMRES-OMS(2)} &\multicolumn{3}{c|}{PCG-OHS(2)} &\multicolumn{2}{c}{GMRES-OMS(2)}  \\
      $\nu$  &it.   &err.    &it.   &err.   &it.   &err.  & cond.=$\lambda_{max}/\lambda_{min}$  &it.   &err.   \\
    \hline
       \cline{1-10}
    &\multicolumn{9}{c}{central jump test} \\
  0.3         &14 &2.8e-7  &5  &1.5e-7  &13 &1.2e-10  &4.50=3.99/0.89   &5  &1.8e-10    \\
  0.4         &15 &4.6e-7  &6  &2.5e-8  &13 &1.8e-10  &4.50=3.99/0.88   &5  &3.3e-10   \\
  0.49        &17 &1.4e-6  &7  &4.4e-7  &15 &8.0e-11  &5.18=3.99/0.77   &6  &1.3e-10    \\
  0.499       &18 &1.1e-6  &8  &5.3e-7  &17 &1.0e-10  &6.86=3.99/0.58   &7  &7.1e-11    \\
  0.4999      &18 &1.2e-6  &8  &7.0e-7  &18 &8.9e-11  &7.49=3.99/0.53   &7  &1.7e-10   \\
  0.49999     &18 &1.2e-6  &8  &7.2e-7  &18 &1.1e-11  &7.55=3.99/0.53   &7  &1.2e-10    \\
    \hline
        \cline{1-10}
    &\multicolumn{9}{c}{checkerboard test} \\
               &18 &1.7e-6  &7  &e-7      &20 &9.4e-11   &8.62=3.99/0.46 &8 &3.6e-10  \\
  \hline
  \end{tabular}
\caption{$Q_n - Q_{n-2}$ spectral elements with $n=2$.  Robustness with respect to jumps of the elliptic coefficients. Saddle point and positive definite reformulations using both OHS(2) and OMS(2) preconditioners. Iteration counts, errors, condition numbers and extreme eigenvalues for an increasing Poisson ratio $\nu \rightarrow \frac{1}{2}$. Fixed $N = 4 \times 4,\ H/h = 4$, overlap $\delta = h$. 
}
\label{jump_test_sem}
\end{center}
\end{table}

\begin{table}[t]
\begin{center}
 \centering
 \begin{tabular}{c|ccc|cc}
 \hline
& \multicolumn{5}{c}{pos. def. reformulation}\\
 & \multicolumn{3}{c}{PCG-OHS(2)}
& \multicolumn{2}{c}{GMRES-OMS(2)}\\
 {$n$} &it. &err. & cond.=$\lambda_{max}/\lambda_{min}$ &it. &err. \\
 \hline
 3 &33 &5.5e-8 &29.89=4.000/1.33e-1 &13 &1.0e-6 \\
 4 &32 &2.1e-7 &29.99=4.000/1.33e-1 &12 &8.2e-6 \\
 5 &31 &4.4e-7 &29.99=4.000/1.33e-1 &12 &4.4e-7 \\
 6 &32 &3.5e-7 &29.99=4.000/1.33e-1 &13 &4.4e-6 \\
 7 &32 &5.0e-7 &29.99=4.000/1.33e-1 &12 &1.3e-5 \\
 8 &32 &6.6e-7 &29.99=4.000/1.33e-1 &12 &6.7e-6 \\
 \hline
 \end{tabular}
\caption{$Q_n-Q_{n-2}$ spectral elements. $n$-independence of PCG-OHS(2) and PCG-OMS(2) for the positive definite reformulation. Iteration counts, errors, condition numbers and extreme eigenvalues for increasing polynomial degree $n$. Fixed $N=3\times 3$, $H/h=4,\ \nu=0.4999$, overlap size $\delta = h$.
}\label{degree-pcg-ohs-oms_sem}
\end{center}
\end{table}

\section*{Acknowledgements}

The authors would like to thank Olof Widlund for his advices and for stimulating discussions.

\end{document}